\documentclass[a4paper,12pt]{article}

\usepackage{amsmath}
\usepackage{ntheorem}
\usepackage{amssymb}
\usepackage{latexsym}
\usepackage{url,mathrsfs}
\usepackage[all]{xy}\UseComputerModernTips
\usepackage{graphicx}

\theoremstyle{plain}
\newtheorem{proposition}{Proposition}[section]
\newtheorem{theorem}[proposition]{Theorem}

\newtheorem{corollary}[proposition]{Corollary}

\theoremstyle{plain}
\theorembodyfont{\normalfont\rm}
\newtheorem{definition}[proposition]{Definition}
\newtheorem{remark}[proposition]{Remark}

\theoremstyle{nonumberplain}
\theoremheaderfont{\normalfont\itshape}
\theoremseparator{.}
\theorembodyfont{\normalfont}

\newcommand{\qed}{\hfill $\Box$}

\makeatletter

\@addtoreset{equation}{section}
\makeatother

\newcommand{\ZZ}{{\mathbb Z}}

\newcommand{\CC}{{\mathbb C}}
\newcommand{\PP}{{\mathbb P}}

\renewcommand{\d}{{\rm dim}}

\newcommand{\e}{\varepsilon}

\newcommand{\id}{{\rm id}}
\newcommand{\supp}{{\rm supp}}

\newcommand{\Int}{{\rm Int}}
\newcommand{\Ker}{{\rm Ker}}

\renewcommand{\Im}{{\rm Im}}

\newcommand{\Gr}{{\rm Gr}}

\newcommand{\Db}{{\bf D}^{b}}
\newcommand{\Dbc}{{\bf D}_{c}^{b}}

\newcommand{\F}{{\cal F}}
\newcommand{\G}{{\cal G}}

\renewcommand{\L}{{\mathcal L}}

\newcommand{\tl}[1]{\widetilde{#1}}

\newcommand{\simto}{\overset{\sim}{\longrightarrow}}

\newcommand{\dsum}{\displaystyle \sum}

\renewcommand{\(}{\left(}
\renewcommand{\)}{\right)}

\newcommand{\longhookrightarrow}{\DOTSB\lhook\joinrel\longrightarrow}

\def\cf{cf.\kern.3em}
\def\eg{e.g.\kern.3em}

\begin{document}

\title{On the sizes of the Jordan blocks of \\monodromies at infinity }

\author{Yutaka \textsc{Matsui}\footnote{Department of Mathematics, Kinki University, 3-4-1, Kowakae, Higashi-Osaka, Osaka, 577-8502, Japan. E-mail: matsui@math.kindai.ac.jp} \and Kiyoshi \textsc{Takeuchi}\footnote{Institute of Mathematics, University  of Tsukuba, 1-1-1, Tennodai, Tsukuba, Ibaraki, 305-8571, Japan. E-mail: takemicro@nifty.com}}

\date{}

\sloppy

\maketitle

\begin{abstract}
We obtain general upper bounds of the sizes and the numbers of Jordan blocks for the eigenvalues $\lambda \not= 1$ in the monodromies at infinity of polynomial maps. 
\end{abstract}

\section{Introduction}\label{sec:1}

In this paper we study the upper bounds of the sizes and the numbers of Jordan blocks in the monodromies at infinity of general polynomial maps. First we recall the definition of monodromies at infinity. After two fundamental papers \cite{Broughton} and \cite{S-T-1}, many authors studied the global behavior of polynomial maps $f \colon \CC^n \longrightarrow \CC$. For a polynomial map $f \colon \CC^n \longrightarrow \CC$, it is well-known that there exists a finite subset $B \subset \CC$ such that the restriction
\begin{equation}
\CC^n \setminus f^{-1}(B) \longrightarrow \CC \setminus B
\end{equation}
of $f$ is a locally trivial fibration. We denote by $B_f$ the smallest subset $B \subset \CC$ satisfying this condition. Let $C_R=\{x\in \CC\ |\ |x|=R\}$ ($R\gg 0$) be a sufficiently large circle in $\CC$ such that $B_f\subset \{x \in \CC\ |\ |x|<R\}$. Then by restricting the locally trivial fibration $\CC^n \setminus f^{-1}(B_f) \longrightarrow \CC \setminus B_f$ to $C_R$ we obtain a geometric monodromy automorphism $\Phi_f^{\infty} \colon f^{-1}(R) \simto f^{-1}(R)$ and the linear maps
\begin{equation}
\Phi_j^{\infty} \colon H^j(f^{-1}(R) ;\CC) \overset{\sim}{\longrightarrow} H^j(f^{-1}(R) ;\CC) \ \ (j=0,1,\ldots)
\end{equation}
associated to it, where the orientation of $C_R$ is taken to be counter-clockwise as usual. We call $\Phi_j^{\infty}$ the (cohomological) monodromies at infinity of $f$. Various formulas for their eigenvalues (i.e. the semisimple parts) were obtained by many authors. In particular, for their expressions in terms of the Newton polyhedra at infinity of $f$, see Libgober-Sperber \cite{L-S} and \cite{M-T-1} etc. Also, some important results on the nilpotent parts of $\Phi_j^{\infty}$ were obtained by Garc{\'i}a-L{\'o}pez-N{\'e}methi \cite{L-N-2} and Dimca-Saito \cite{D-S-1} etc. For example, Dimca-Saito \cite{D-S-1} obtained an upper bound of the sizes of Jordan blocks for the eigenvalue $1$ in $\Phi_j^{\infty}$. Recently in \cite{M-T-2} we obtained very explicit formulas which express the Jordan normal forms of $\Phi_j^{\infty}$ in terms of the Newton polyhedra at infinity of $f$ (see \cite{M-T-3} and \cite{E-T} for the further developments). However they are applicable only to convenient polynomials $f$ which are non-degenerate at infinity. By a result of Broughton \cite{Broughton}, such polynomials are tame at infinity in the sense of Kushnirenko \cite{Kushnirenko}. In this paper, without assuming that $f$ is tame at infinity, we obtain a general upper bound of the sizes of Jordan blocks for each eigenvalue $\lambda \not= 1$ in $\Phi_j^{\infty}$, which is similar to the one for the eigenvalue $1$ in Dimca-Saito \cite{D-S-1}. Moreover we also give an upper bound of the numbers of such Jordan blocks with the maximal possible size $j+1$ in $\Phi_j^{\infty}$. In the course of our proof, the methods in their another paper \cite{D-S-2} will be effectively used.

\section{Monodromies at infinity}\label{sec:2}

In this section, we recall some basic definitions on monodromies at infinity. Let $f(x) \in \CC [x_1, x_2, \ldots, x_n]$ be a polynomial on $\CC^n$. Then as we explained in Introduction, there exist a locally trivial fibration $\CC^n \setminus f^{-1}(B_f) \longrightarrow \CC \setminus B_f$ and the linear maps
\begin{equation}
\Phi_j^{\infty} \colon H^j(f^{-1}(R) ;\CC) \overset{\sim}{\longrightarrow} H^j(f^{-1}(R) ;\CC) \ \ (j=0,1,\ldots)
\end{equation}
($R \gg 0$) associated to it. To study the monodromies at infinity $\Phi_j^{\infty}$, we often impose the following natural condition. 

\begin{definition}[\cite{Kushnirenko}]\label{dfn:tame}
Let $\partial f\colon \CC^n \longrightarrow \CC^n$ be the map defined by $\partial f(x)=(\partial_1f(x), \ldots, \partial_n f(x))$. Then we say that $f$ is tame at infinity if the restriction $(\partial f)^{-1}(B(0;\e )) \longrightarrow B(0;\e )$ of $\partial f$ to a sufficiently small ball $B(0;\e )$ centered at the origin $0 \in \CC^n$ is proper.
\end{definition}

The following result is fundamental in the study of monodromies at infinity.

\begin{theorem}[Broughton \cite{Broughton} and Siersma-Tib{\u a}r \cite{S-T-1}]\label{tame}
Assume that $f$ is tame at infinity. Then the generic fiber $f^{-1}(c)$ ($c \in \CC \setminus B_f$) has the homotopy type of the bouquet of $(n-1)$-spheres. In particular, we have
\begin{equation}
H^j(f^{-1}(c);\CC)=0 \quad (j \neq 0, n-1).
\end{equation}
\end{theorem}

By this theorem if $f$ is tame at infinity, $\Phi_{n-1}^{\infty}$ is the only non-trivial monodromy at infinity. Many authors studied tame polynomials. However, in this paper we do not assume the tameness at infinity of $f$ and study the general properties of the monodromies at infinity $\Phi_j^{\infty}$.

The following general result is often called the monodromy theorem.

\begin{theorem}
For $\lambda \in \CC \setminus \{1\}$ the sizes of Jordan blocks for the eigenvalue $\lambda$ in $\Phi_j^{\infty}$ are $\leq j+1$. 
\end{theorem}

\section{Some properties of the nearby cycle functor}\label{sec:3}

The nearby cycle functor introduced by Deligne will play an important role in this paper. In this paper, we essentially follow the terminology in \cite{Dimca} and \cite{K-S}. For example, for an algebraic variety $X$ over $\CC$, we denote by $\Db(X)$ the derived category of bounded complexes of sheaves of $\CC_X$-modules on $X$, by $\Dbc(X)$ the full subcategory of $\Db(X)$ consisting of bounded complexes of sheaves whose cohomology sheaves are constructible and by ${\rm Perv}(X)$ the category of perverse sheaves on $X$. For the detail, see \cite{Dimca}, \cite{H-T-T}, \cite{K-S}, \cite{Saito-1} and \cite{Saito-2}.

\begin{definition}
Let $X$ be an algebraic variety over $\CC$ and $f \colon X \longrightarrow \CC$ a non-constant regular function on $X$. Set $X_0:= \{x\in X\ |\ f(x)=0\} \subset X$ and let $i_X \colon X_0 \longhookrightarrow X$, $j_X \colon X \setminus X_0 \longhookrightarrow X$ be inclusions. Let $p \colon \tl{\CC^*} \longrightarrow \CC^*$ be the universal covering of $\CC^* =\CC \setminus \{0\}$ ($\tl{\CC^*} \simeq \CC$) and consider the Cartesian square
\begin{equation}\label{eq:pre-1}
\xymatrix@R=2.5mm@C=2.5mm{
\tl{X \setminus X_0} \ar[rr] \ar[dd]^{p_X} & &\tl{\CC^*} \ar[dd]^p \\
 & \Box & \\
X \setminus X_0 \ar[rr]^f & & \CC^*.}
\end{equation}
Then for $\G \in \Db(X)$ we set
\begin{equation}
\psi_f(\G) := i_X^{-1}R(j_X \circ p_X)_*(j_X \circ p_X)^{-1}\G \in \Db(X_0)
\end{equation}
and call it the nearby cycle of $\G$. 

Let us denote by ${\rm Deck}(\tl{\CC^*},\CC^*)\simeq \ZZ$ the group of deck transformations of the covering map $p \colon \tl{\CC^*}\longrightarrow \CC^*$. The action of a generator $1\in \ZZ$ of ${\rm Deck}(\tl{\CC^*},\CC^*)\simeq \ZZ$ on $\tl{X\setminus X_0}$ induces an automorphism $\Phi(\G)$ of $\psi_f(\G)$
\begin{equation}
\Phi(\G) \colon \psi_f(\G)\simto \psi_f(\G).
\end{equation}
We call it the monodromy automorphism of $\psi_f(\G)$.
\end{definition}

Since the nearby cycle functor $\psi_f$ preserves the constructibility, we obtain the functor
\begin{equation}
\psi_f\colon \Dbc(X) \longrightarrow \Dbc(X_0).
\end{equation}
Moreover, since $\psi_f[-1]$ preserves the perversity, we obtain the functor
\begin{equation}
\psi_f[-1] \colon {\rm Perv}(X) \longrightarrow {\rm Perv}(X_0).
\end{equation}

The nearby cycle functor $\psi_f$ generalizes the classical notion of Milnor fibers. Suppose that $X$ is a subvariety of $\CC^m$ and $f \colon X\longrightarrow \CC$ is a non-constant regular function. Then for $x\in X_0$ we can define the local Milnor fiber $F_x$ of $f$ at $x$. We have the following fundamental result (for example see \cite[Proposition 4.2.2]{Dimca}).

\begin{theorem}\label{prp:3-2} For any $\G\in \Dbc(X)$, $x\in X_0$ and $j\in \ZZ$, there exists a natural isomorphism
\begin{equation}
H^j(F_x;\G)\simeq H^j(\psi_f(\G))_x.
\end{equation}
\end{theorem}

Let us recall briefly some results in \cite[Section 1.4]{D-S-2}. Let $X$ be an $n$-dimensional smooth algebraic variety and $f\colon X\longrightarrow \CC$ a non-constant regular function on $X$. Note that $\CC_X[n]$ is a perverse sheaf on $X$ and the mixed Hodge module corresponding to $\CC_X[n]$ is pure of weight $n$. Set $\G:=\CC_X[n-1]$ and $\F:=\psi_f(\G)\in {\rm Perv}(X_0)$. The monodromy automorphism $\Phi:=\Phi(\G)$ induces the following canonical decomposition
\begin{equation}
\F=\bigoplus_{\lambda\in \CC}\F_{\lambda},
\end{equation}
where we set
\begin{equation}
\F_{\lambda}:=\Ker \left[ (\lambda \cdot \id -\Phi)^N \colon \F \longrightarrow \F \right]\in {\rm Perv}(X_0)
\end{equation}
for $N \gg 0$. Note that for $x\in X_0$ the stalk $H^{j-n+1}(\F_{\lambda})_x$ is isomorphic to the generalized $\lambda$-eigenspace of the classical Milnor monodromy automorphism $H^{j}(F_x;\CC)\simto H^{j}(F_x;\CC)$ by Theorem \ref{prp:3-2}. Let $\Phi |_{\F_{\lambda}}=(\lambda \cdot \id) \times \Phi_u$ be the Jordan decomposition of $\Phi |_{\F_{\lambda}} \colon \F_{\lambda} \longrightarrow \F_{\lambda}$ ($\Phi_u$ is the unipotent part) and set
\begin{equation}
N_{\lambda}:= \dfrac{1}{2\pi \sqrt{-1}}\log{\Phi_u}= \dfrac{1}{2\pi \sqrt{-1}} \dsum_{i=1}^N\dfrac{(-1)^{i+1}}{i} (\Phi_u-\id)^i
\end{equation}
for $N \gg 0$. Then $N_{\lambda}$ is a nilpotent endomorphism of $\F_{\lambda}$. Considering the mixed Hodge module associated with the perverse sheaf $\F_{\lambda}\oplus \F_{\overline{\lambda}}$, the monodromy filtration induced by $N_{\lambda}$ gives the weight filtration $W$ of $\F_{\lambda}$. Recall that $N_{\lambda}$ is strict with respect to the filtration $W$ for the shift $-2$ (see e.g. \cite[Section 1.4]{D-S-1} etc. for the details). Since the mixed Hodge module corresponding to $\CC_X[n]$ is pure of weight $n$, we have the following isomorphism
\begin{equation}\label{eq:3-1}
N_{\lambda}^i \colon \Gr^W_{n-1+i}(\F_{\lambda})\simto \Gr^W_{n-1-i}(\F_{\lambda})
\end{equation}
for any $i\geq 0$ (\cite[\S 5]{Saito-1}). Let us define the primitive part $P\Gr_{n-1+i}^W(\F_{\lambda})$ by
\begin{equation}
P\Gr_{n-1+i}^W(\F_{\lambda})
:=\begin{cases}\Ker[ N_{\lambda}^{i+1}\colon \Gr_{n-1+i}^W(\F_{\lambda})\longrightarrow \Gr_{n-3-i}^W(\F_{\lambda})] & (i\geq 0),\\ 0& (i<0).
\end{cases}
\end{equation}
Then by \eqref{eq:3-1} for each $k$ we have the primitive decomposition of $\Gr_k^W(\F_{\lambda})$:
\begin{equation}\label{PD}
\Gr_k^W(\F_{\lambda})=\bigoplus_{i\geq 0}N_{\lambda}^i\(P\Gr_{k+2i}^W(\F_{\lambda})\).
\end{equation}

In this paper, we will use the following geometric description of the primitive part $P\Gr_{n-1+i}^W(\F_{\lambda})$ in \cite[3.3]{Saito-2}. From now on, let us assume that $X_0=f^{-1}(0)$ is a strictly normal crossing divisor in $X$. Namely, we assume that $X_0$ is a normal crossing divisor whose irreducible components $D_1,\ldots, D_m$ are smooth. For $1\leq i\leq m$, let $a_i>0$ be the order of the zeros of $f$ along $D_i$. For $\lambda \in \CC$, we set $R_{\lambda}:=\{1\leq i\leq m\ |\ \lambda^{a_i}=1\}\subset \{1,\ldots,m\}$. Moreover, for a non-empty subset $I\subset R_{\lambda}$ we set 
\begin{equation}
D_I:=\bigcap_{i\in I}D_i, \hspace{10mm}U_I:=D_I\setminus\( \bigcup_{i\notin R_{\lambda}}D_i\).
\end{equation}
For a non-empty subset $I\subset R_{\lambda}$, let $\L_{\lambda, I}$ be a local system of rank $1$ on $U_I$ whose monodromy around the divisor $D_i$ for $i\notin R_{\lambda}$ is defined by the multiplication by $\lambda^{-a_i}(\neq 1)$. Then we have the following decomposition of the primitive part $P\Gr_{n-1+i}^W(\F_{\lambda})$ 
\begin{equation}\label{eq:pgr}
P\Gr_{n-1+i}^W(\F_{\lambda})\simeq\bigoplus_{\begin{subarray}{c}I\subset R_{\lambda}\\ \sharp I=i+1\end{subarray}} (j_I)_!\L_{\lambda, I}[n-1-i],
\end{equation}
where $j_I \colon U_I \longhookrightarrow X_0$ is the natural inclusion. Note that we have an isomorphism
\begin{equation}
(j_I)_!\L_{\lambda, I}[n-1-i]\simeq R(j_I)_*\L_{\lambda,I}[n-1-i].
\end{equation}
By \eqref{eq:pgr}, for $i\geq \max\{ \sharp I\ |\ I\subset R_{\lambda}, \ D_I\neq \emptyset\}$, we have $\Gr_{n-1+i}^W(\F_{\lambda})=0$ and $N_{\lambda}^i=0$.

\section{Main results}\label{sec:4}

In this section, without assuming that $f$ is tame at infinity, we prove some general results on the sizes and the numbers of the Jordan blocks in the monodromies at infinity $\Phi_j^{\infty}$ of $f$. Let $X$ be a smooth compactification of $\CC^n$. Then by eliminating the points of indeterminacy of the meromorphic extension of $f$ to $X$ we obtain a commutative diagram
\begin{equation}
\xymatrix{
\CC^n \ar@{^{(}->}[r]^{\iota} \ar[d]_f & \tl{X} \ar[d]^g\\
\CC \ar@{^{(}->}[r]^j & \PP^1}
\end{equation}
such that $g$ is a proper holomorphic map and $\tl{X} \setminus \CC^n$, $Y:=g^{-1}( \infty ) \subset \tl{X} \setminus \CC^n$ are strict normal crossing divisors in $\tl{X}$. See e.g. Sabbah \cite{Sabbah} and \cite[Section 4]{M-T-2} etc. Let us define an open subset $\Omega$ of $\tl{X}$ by
\begin{equation}
\Omega=\Int (\iota(\CC^n) \sqcup Y)
\end{equation}
and set $U=\Omega \cap Y$. Then $U$ (resp. the complement of $\Omega$ in $\tl{X}$) is a normal crossing divisor in $\Omega$ (resp. $\tl{X}$). In this situation, the main result of Dimca-Saito \cite{D-S-1} can be stated as follows.

\begin{theorem}{\rm 
\bf(\cite[Theorem 0.1]{D-S-1})}
Let $F_1, F_2, \ldots, F_l$ be the irreducible components of $\tl{X} \setminus \CC^n$ contained in $\tl{X} \setminus \Omega$. Assume that for generic complex numbers $c \in \CC$ the closures $\overline{f^{-1}(c)}$ of $f^{-1}(c)$ in $\tl{X}$ are smooth and intersect $\bigcap_{i \in I}F_i$ for any $I \subset \{ 1,2, \ldots, l \}$ transversally. By taking such a complex number $c \in \CC$ we set
\begin{equation}
K= \max_{p \in (\tl{X} \setminus \Omega) \cap \overline{f^{-1}(c)}} \(\sharp \{ F_i \ |\ p \in F_i\}\).
\end{equation}
Then the size of the Jordan blocks for the eigenvalue $1$ of the monodromies at infinity $\Phi_j^{\infty} \colon H^j(f^{-1}(R);\CC)\simto H^j(f^{-1}(R);\CC)$ ($R\gg 0$, $j=0,1,\ldots$) is bounded by $K$.
\end{theorem}

By using Saito's mixed Hodge modules in a different way, we can prove a similar result also for the eigenvalues $\lambda \in \CC \setminus \{1\}$ of $\Phi_j^{\infty}$ as follows. Recall that the size of the Jordan blocks for such eigenvalues in $\Phi_j^{\infty}$ is bounded by $j+1$ by the monodromy theorem. Let $E_1, E_2, \ldots, E_k$ be the irreducible components of the normal crossing divisor $U=\Omega \cap Y$ in $\Omega \subset \tl{X}$ and for $1 \leq i \leq k$ let $b_i>0$ be the order of the poles of $f$ along $E_i$. For a subset $I \subset \{1,\ldots, k \}$ we set $E_I=\bigcap_{i \in I} E_i$. Moreover for $\lambda \in \CC$ we set
\begin{equation}
R_{\lambda}=\{ 1\leq i \leq k \ |\ \lambda^{b_i}=1 \} \subset \{1,\ldots, k \}.
\end{equation}

\begin{theorem}\label{thm:6-2}
Assume that $\lambda \in \CC \setminus \{1\}$.
\begin{enumerate}
\item  We set
\begin{equation}
K_{\lambda}= \max_{p \in U} \(\sharp \{E_i \ |\ p\in E_i \quad \text{and} \quad \lambda^{b_i}=1\}\).
\end{equation}
Then for any $0 \leq j \leq n-1$ the size of the Jordan blocks for the eigenvalue $\lambda$ in $\Phi_j^{\infty}$ is bounded by $K_{\lambda}$.
\item For $0 \leq j \leq n-1$, we set
\begin{equation}
S(\lambda)_j= \{ I \subset R_{\lambda} \ |\ \sharp I =j+1 \quad \text{and} \quad E_I \neq \emptyset \}.
\end{equation}
Moreover for each $I \in S(\lambda)_j$ denote by $c_I$ the number of the connected components of $E_I$ which do not intersect $E_i$ for any $i \notin R_{\lambda}$. Then the number of the Jordan blocks for the eigenvalue $\lambda$ with the maximal possible size $j+1$ in $\Phi_j^{\infty}$ is bounded by $\sum_{I \in S(\lambda)_j} c_I$.
\end{enumerate}
\end{theorem}

\begin{corollary}
We set
\begin{equation}
S(\lambda)=\left\{ I \subset R_{\lambda} \ |\ \sharp I =n \quad \text{and} \quad E_I \neq \emptyset\right\}
\end{equation}
and for each $I \in S(\lambda)$ denote by $n_I$ the cardinality of the discrete (hence finite) set $E_I$. Then the number of the Jordan blocks for $\lambda \in \CC\setminus \{1\}$ with the maximal possible size $n$ in $\Phi_{n-1}^{\infty}$ is bounded by $\sum_{I \in S(\lambda)} n_I$.
\end{corollary}

\noindent{\bf Proof of Theorem \ref{thm:6-2}}
 \ Set $\tl{g} =\frac{1}{f}$. Then for $R \gg 0$ we can easily prove the isomorphisms
\begin{equation}
H^j(f^{-1}(R);\CC) \simeq H^j(Y; \psi_{\tl{g}}(R\iota_*\CC_{\CC^n}))\simeq H^j(U; \psi_{\tl{g}}(\CC_{\tl{X}})).
\end{equation}
Now let us consider the nearby cycle perverse sheaf $\F =\psi_{\tl{g}}(\CC_{\tl{X}}[n-1])\in \Dbc(Y)$ on the normal crossing divisor $Y$ and its monodromy automorphism
\begin{equation}
\Phi :=\Phi(\CC_{\tl{X}}[n-1] ) \colon \F \simto \F.
\end{equation}
Then for $R \gg 0$ we have a commutative diagram
\begin{equation}
\xymatrix@R=5mm@C=15mm{
H^j(f^{-1}(R);\CC)\ar[r]^{\Phi_j^{\infty}} \ar@{-}[d]^{\wr} & H^j(f^{-1}(R);\CC) \ar@{-}[d]^{\wr}\\
H^{j-n+1}(U;\F) \ar[r]^{\Phi} & H^{j-n+1}(U;\F).}
\end{equation}
Moreover there exists a canonical decomposition
\begin{equation}
\F=\bigoplus_{\lambda \in \CC} \F_{\lambda},
\end{equation}
where we set $\F_{\lambda}=\Ker \left[ (\lambda \cdot \id -\Phi)^N \colon \F \longrightarrow \F \right]$ for $N \gg 0$. Therefore, for the given $\lambda \in \CC \setminus\{1\}$ the generalized eigenspace for the eigenvalue $\lambda$ in $\Phi_j^{\infty} \colon H^j(f^{-1}(R);\CC) \simto H^j(f^{-1}(R);\CC)$ ($R \gg 0$) is isomorphic to $H^{j-n+1}(U;\F_{\lambda})$. Now let $\Phi |_{\F_{\lambda}}=(\lambda \cdot \id) \times \Phi_u$ be the Jordan decomposition of $\Phi |_{\F_{\lambda}} \colon \F_{\lambda} \longrightarrow \F_{\lambda}$ ($\Phi_u$ is the unipotent part) and set
\begin{equation}
N_{\lambda}= \dfrac{1}{2\pi \sqrt{-1}}\log{\Phi_u}
= \dfrac{1}{2\pi \sqrt{-1}} \dsum_{i=1}^N\dfrac{(-1)^{i+1}}{i} 
(\Phi_u-\id)^i
\end{equation}
for $N \gg 0$. Then $N_{\lambda}$ is a nilpotent endomorphism of the perverse sheaf $\F_{\lambda}$ and there exists an automorphism $M_{\lambda} \colon \F_{\lambda} \longrightarrow \F_{\lambda}$ such that
\begin{equation}
\Phi_u -\id =N_{\lambda}M_{\lambda}=M_{\lambda}N_{\lambda}.
\end{equation}
This implies that if $(N_{\lambda})^i=0$ for some $i \geq 1$ then $(\lambda \cdot\id -\Phi |_{\F_{\lambda}})^i=\lambda^i(\id -\Phi_u)^i=0$ and the size of the Jordan blocks for the eigenvalue $\lambda$ in $\Phi_j^{\infty} \colon H^j(f^{-1}(R);\CC) \simto H^j(f^{-1}(R);\CC)$ is $\leq i$. Let $W$ be the weight filtration of the mixed Hodge module associated with the perverse sheaf $\F_{\lambda}\oplus \F_{\overline{\lambda}}$. Then the assertion (i) follows from the geometric description of the primitive decomposition \eqref{PD} of the graded module $\Gr^W(\F_{\lambda})$ in Section \ref{sec:3}. Finally let us prove (ii). By the above argument, the number of the Jordan blocks for the eigenvalue $\lambda$ with the maximal possible size $j+1$ in $\Phi_j^{\infty}$ is equal to
\begin{equation}\label{equ}
\dim \(\Im \left[ H^{j-n+1}(U;\F_{\lambda}) \overset{N_{\lambda}^j}{\longrightarrow} H^{j-n+1}(U;\F_{\lambda}) \right] \).
\end{equation}
Let $\G$ be the subobject $\Im N_{\lambda}^j$ of $\F_{\lambda}$ in the category of perverse sheaves on $Y$. Then we have a commutative diagram
\begin{equation}
\xymatrix@R=5mm@C=5mm{
H^{j-n+1}(U;\F_{\lambda}) \ar[rd] \ar[rr]^{N_{\lambda}^j} & & 
H^{j-n+1}(U;\F_{\lambda})\\
&H^{j-n+1}(U;\G ), \ar[ur] &}
\end{equation}
and the number \eqref{equ} is bounded by $\d H^{j-n+1}(U;\G )$. Let us set $\G^{\prime}=\G \cap W_{n-j-1}\F_{\lambda}$ and $\G^{\prime \prime}=\G \cap W_{n-j-2}\F_{\lambda}$. Then by the structure of the primitive decomposition \eqref{PD} of $\Gr^W(\F_{\lambda})$ we have $\d \ \supp (\G/\G^{\prime}) \leq n-j-2$ and $\d \ \supp \G^{\prime \prime} \leq n-j-2$. Here we used the strictness of $N_{\lambda}$ with respect to the filtration $W$ for the shift $-2$. Hence we obtain
\begin{equation}
H^i(U; \G/\G^{\prime})=H^i(U; \G^{\prime \prime}) =0 \quad \text{for any} \ i < j-n+2. 
\end{equation}
Then by the exact sequence of perverse sheaves
\begin{equation}
0 \longrightarrow \G^{\prime} \longrightarrow \G \longrightarrow \G/\G^{\prime} \longrightarrow 0
\end{equation}
we obtain an isomorphism
\begin{equation}
H^{j-n+1}(U;\G^{\prime} ) \simeq H^{j-n+1}(U;\G ).
\end{equation}
Moreover it follows from the exact sequence of perverse sheaves
\begin{equation}
0 \longrightarrow \G^{\prime \prime} \longrightarrow \G^{\prime} \longrightarrow \G^{\prime}/\G^{\prime \prime} \simeq \Gr^W_{n-j-1}\F_{\lambda} \longrightarrow 0
\end{equation}
that we have
\begin{equation}
\d H^{j-n+1}(U;\G^{\prime}) \leq \d H^{j-n+1}(U; \Gr^W_{n-j-1}\F_{\lambda}).
\end{equation}
For $I \subset R_{\lambda}$ such that $\sharp I=j+1$, set $U_I =E_I \setminus \( \bigcup_{i \notin R_{\lambda}}E_i \)$. Then by the geometric description of the primitive decomposition \eqref{PD} of $\Gr^W(\F_{\lambda})$ in Section \ref{sec:3} we can easily see that
\begin{equation}
\d H^{j-n+1}(U; \Gr^W_{n-j-1}\F_{\lambda}) \leq \d \(\bigoplus_{\begin{subarray}{c}I \subset R_{\lambda} \\ 
\sharp I=j+1\end{subarray}} \varGamma(U_I; \L_{\lambda,I})\),
\end{equation}
where $\L_{\lambda,I}$ is a local system of rank one on $U_I$ whose monodromy around the divisor $E_i$ ($i \notin R_{\lambda}$) is given by the multiplication by $\lambda^{-b_i} (\neq 1)$. If a connected component $E_{I,r}$ of $E_I$ intersects $E_i$ for some $i \not\in R_{\lambda}$ we have $\varGamma(U_I\cap E_{I,r};\L_{\lambda,I})=0$. Therefore the assertion (ii) follows. This completes the proof. \qed

\begin{remark}
By Theorem \ref{thm:6-2} (ii) it seems that for $0 \leq j \leq n-2$ and $\lambda \in \CC \setminus \{1\}$ there is no Jordan block for the eigenvalue $\lambda$ with the maximal possible size $j+1$ in $\Phi_j^{\infty}$ in general. This implies that the generalized $\lambda$-eigenspaces of the monodromies at infinity $\Phi_j^{\infty}$ ($0 \leq j \leq n-2$) are much simpler than that of the top one $\Phi_{n-1}^{\infty}$. For similar results in the case of local Milnor monodromies, see Dimca-Saito \cite{D-S-2}. 
\end{remark}

\end{document}